\documentclass[12pt]{article}
\usepackage{amsfonts}
\usepackage{amsmath,amssymb,amsthm}
\usepackage[usenames]{color}
\usepackage{mathrsfs}
\usepackage{syntonly}
\usepackage[top=1in,bottom=1.4in,left=1in,right=1in]{geometry}
\usepackage[all,cmtip]{xy}
\usepackage{bbm}
\usepackage{mathrsfs}
\usepackage{dsfont}
\usepackage{enumitem}
\usepackage{fouridx}
\usepackage[stable]{footmisc}
\usepackage{array}
\usepackage{caption}
\usepackage{latexsym}
\usepackage{extarrows}
\usepackage{amsfonts}
\usepackage{amsmath}
\usepackage{amsthm}
\usepackage{amssymb}
\usepackage{latexsym}
\usepackage{indentfirst}
\allowdisplaybreaks

\newcommand{\C}{\mathbb{C}}
\newcommand{\Z}{\mathbb{Z}}
\newcommand{\N}{\mathbb{N}}

\newcommand{\g}{\mathfrak{g}}

\newcommand{\yn}{Y(\mathfrak{sl}_{l+1})}
\newcommand{\yg}{Y(\mathfrak{g})}
\newcommand{\ysl}{Y(\mathfrak{sl}_2)}

\newcommand{\nyn}{\mathfrak{sl}_{l+1}}
\newcommand{\nyo}{\mathfrak{so}(2l,\C)}

\newcommand{\nysp}{\mathfrak{sp}(2l,\C)}

\newcommand{\nyso}{\mathfrak{so}(2l+1,\C)}

\numberwithin{equation}{section}
\newtheorem{theorem}{Theorem}[section] 
\newtheorem{proposition}[theorem]{Proposition}
\newtheorem{corollary}[theorem]{Corollary}
\newtheorem{lemma}[theorem]{Lemma}
\newtheorem{definition}[theorem]{Definition}
\newtheorem{remark}[theorem]{Remark}

\begin{document}


\title{\Large\bf  Local Weyl modules and cyclicity of tensor products for  Yangians}
\author{{Yilan Tan, Nicolas Guay}}

\date{} 

\maketitle

\vspace{17 mm}

\begin{abstract}
We provide a sufficient condition for the cyclicity of an ordered tensor product $L=V_{a_1}(\omega_{b_1})\otimes V_{a_2}(\omega_{b_2})\otimes\ldots\otimes V_{a_k}(\omega_{b_k})$ of fundamental representations of the Yangian $\yg$. When $\g$ is a classical simple Lie algebra, we make the cyclicity condition concrete, which leads to an irreducibility criterion for the ordered tensor product $L$. In the case when $\g=\nyn$, a sufficient and necessary condition for the irreducibility of the ordered tensor product $L$ is obtained. The cyclicity of the ordered tensor product $L$ is closely related to the structure of the local Weyl modules of $\yg$. We show that every local Weyl module is isomorphic to an ordered tensor product of fundamental representations of $\yg$.
\end{abstract}





\newpage
\section{Introduction}
Yangians and quantum affine algebras form two of the most important classes of quantum groups. The representation theory of Yangians has many applications in mathematics and physics. For instance, from any finite-dimensional representation of a Yangian, one can construct a $R$-matrix which is a  solution of the quantum Yang-Baxter equation \cite{ChPr4}. As stated in \cite{NT}, ``the physical data such as mass formula, fusion angle, and the spins of integrals of motion can be extracted from the Yangian highest weight representations;" see also \cite{ChPr8}.

One fundamental problem regarding the highest weight representation of a Yangian $\yg$ is to determine the structure of the finite-dimensional irreducible representations.
Unfortunately, this is unknown in most cases except for $\g=\mathfrak{sl}_{l+1}$ where a lot is known thanks to the work of S. Khoroshkin, M. Nazarov et al., see for instance \cite{KhNaPa} and related papers. In this article, we focus on the cases when $\g$ is a finite-dimensional simple Lie algebra over $\C$ of rank $l$.
It was proved in \cite{Dr2} that every finite-dimensional irreducible representation $L$ of $\yg$ is a highest weight representation, and its highest weight is encapsulated by an $l$-tuple of monic polynomials $\pi=\big(\pi_1(u),\ldots,\pi_l(u)\big)$, where $\pi_i(u)=\prod\limits_{j=1}^{m_i}(u-a_{i,j})$; $\pi_i(u)$ is called a Drinfeld polynomial.
An irreducible representation is called fundamental if there is an $i$ such that $\pi_i(u)=u-a$ and $\pi_{j}(u)=1$ for any $j\neq i$. In this case, we denote the representation by $V_{a}(\omega_i)$. 
V. Chari and A. Pressley showed in \cite{ChPr4} that the finite-dimensional irreducible representation of $\yg$ associated to $\pi$ is a subquotient of a tensor product of fundamental representations
$\bigotimes\limits_{i,j} V_{a_{i,j}}(\omega_i)$,
where the tensor factors take any order. 

The finite-dimensional representation theory of the quantum affine algebra $U_q(\hat{\g})$ over $\C(q)$ is an analogue of the one for $\yg$, where $q$ is an indeterminate and $\C(q)$ is the field of rational functions in $q$ with complex coefficients. Finite-dimensional irreducible representations of the quantum affine algebra $U_q(\hat{\g})$  are parameterized by $l$-tuples of polynomials $\operatorname{P}=\big(P_1(u),\ldots, P_l(u)\big)$, where $P_i(0)=1$.  An irreducible representation is called fundamental if there is an $i$ such that $\pi_i(u)=1-cu$ for $c\neq 0$ and $\pi_{j}(u)=1$ for any $j\neq i$. In this case, we denote the representation by $L_{c}(\omega_i)$. It is well known that the finite-dimensional irreducible representation associated to $\operatorname{P}$ is a subquotient of $\bigotimes\limits_{ij}L_{c_{i,j}}(\omega_{i})$, where the tensor factors take any order and the $c_{i,j}$ are roots of $P_i(u)$. In \cite{AkKa}, the authors conjectured that if for any ${c_{i}}$ and ${c_{j}}$, $1\leq i\neq j\leq k$, $\frac{c_{i}}{c_{j}}$ does not have a pole at $q=0$, then the ordered tensor product $\widetilde{L}=L_{c_1}(\omega_{d_1})\otimes L_{c_2}(\omega_{d_2})\otimes \ldots\otimes L_{c_k}(\omega_{d_k})$ is irreducible, and proved this conjecture in the case of type $A_n^{(1)}$ and $C_n^{(1)}$. This conjecture was also proved by E. Frenkel and E. Mukhin \cite{FrMu} using the $q$-characters method, by M. Varagnolo and E. Vasserot \cite{VaVa} via quiver varieties when $\g$ is simply-laced, and by M. Kashiwara in \cite{Ka} through crystal bases. This result was generalized in \cite{Ch3} by V. Chari who gave a sufficient condition for the cyclicity of the tensor product of Kirillov-Reshetikhin modules, i.e., finite-dimensional irreducible representations associated to an $l$-tuple of polynomials $\operatorname{P}$ such that the roots of $P_i(u)$ form a `$q$-string' and $P_j(u)=1$, for all $j\neq i$. This condition is obtained by considering a braid group action on the imaginary root vectors. However, there is no braid group action available for the representation of $\yg$; there is no analogue of the sixth relation in (2.4) in \cite{Ch3} for Yangians.

In order to have a better understanding of the category of finite-dimensional irreducible representations of quantum affine algebras associated to an $l$-tuple of polynomials $\operatorname{P}$, the local Weyl module $W(\operatorname{P})$ was introduced in \cite{ChPr1}. The module $W(\operatorname{P})$ has a nice universal property:  any finite-dimensional highest weight representation of $U_q(\hat{\g})$ associated to $\operatorname{P}$ is a quotient of $W(\operatorname{P})$. It is known that $W(\operatorname{P})$ is isomorphic to an ordered tensor product of fundamental representations of $U_q(\hat{\g})$. A proof of this fact can be found in \cite{ChMo2}. The notion of a local Weyl module has been extended to the finite-dimensional representations of current algebras \cite{ChLo, FoLi, Na}, twisted loop algebras \cite{ChFoSe}, and current Lie algebras on affine varieties \cite{FeLo1}.

In this paper, we study the local Weyl modules $W(\pi)$ and the cyclicity condition for a tensor product of fundamental representations of $\yg$. In Section 2, we introduce the Yangian $\yg$ and its finite-dimensional representations. In Section 3, we give the definition of the local Weyl module $W(\pi)$ of $\yg$ and prove that its dimension is bounded by the dimension of the local Weyl module $W(\lambda)$ of the current algebra $\g[t]$, where $\lambda=\sum\limits_{i\in I} m_i\omega_i$. In Section 4, we provide the structure of the local Weyl module $W(\pi)$ of $\ysl$ by showing that $W(\pi)$ is isomorphic to an ordered tensor product of fundamental representations of $\ysl$, see Theorem \ref{wmfsl2}. In Section 5, we first provide a cyclicity condition for an ordered tensor product of fundamental representations $L=V_{a_1}(\omega_{b_1})\otimes V_{a_2}(\omega_{b_2})\otimes\ldots\otimes V_{a_k}(\omega_{b_k})$, see Theorem \ref{S5T1Liahw}.  We next make the cyclicity condition concrete when $\g$ is classical, see Theorem \ref{mtoyspl}. Recall that a finite-dimensional representation $L$ of $\yg$ is irreducible if both $L$ and the left dual$\ ^t L$ are highest weight representations. Since the left dual of a tensor product of fundamental representations is again a tensor product of fundamental ones, our cyclicity condition for $L$ leads to an irreducibility criterion on $L$, see Theorem \ref{icotpofr}. In particular, when $\g=\nyn$, a sufficient and necessary condition on the irreducibility of $L$ is obtained (See \cite{Mo2} and \cite{NaTa} for more general irreducibility results for $Y(\mathfrak{gl}_{l+1})$). In Theorem \ref{wmiatpsp24}, we show that $W(\pi)$ is isomorphic to an ordered tensor product of fundamental representations. This parallels the results obtained on the structure of the local Weyl modules of quantum affine algebras.
\section{Yangians and their representations}
In this section, we recall the definition of the Yangian $\yg$ and some results concerning its finite-dimensional representations.
\begin{definition}
Let $\g$ be a simple Lie algebra over $\C$ with rank $l$ and let $A=\left(a_{ij}\right)_{i,j\in I}$, where $I=\{1,2,\ldots, l\}$, be its Cartan matrix. Let $D=\operatorname{diag}\left(d_1,\ldots, d_{l}\right)$, $d_i\in \N$, such that $d_1, d_2,\ldots, d_l$ are co-prime and $DA$ is symmetric. The Yangian $Y\left(\mathfrak{g}\right)$ is defined to be the
associative algebra with generators $x_{i,r}^{{}\pm{}}$,
$h_{i,r}$, $i\in I$, $r\in\Z_{\geq 0}$, and the following
defining relations:
\begin{equation*}\label{}
[h_{i,r},h_{j,s}]=0, \qquad [h_{i,0}, x_{j,s}^{\pm}]={}\pm
d_ia_{ij}x_{j,s}^{\pm}, \qquad [x_{i,r}^+, x_{j,s}^-]=\delta_{i,j}h_{i,r+s},
\end{equation*}
\begin{equation*}\label{}
[h_{i,r+1}, x_{j,s}^{\pm}]-[h_{i,r}, x_{j,s+1}^{\pm}]=
\pm\frac{1}{2}d_i
a_{ij}\left(h_{i,r}x_{j,s}^{\pm}+x_{j,s}^{\pm}h_{i,r}\right),
\end{equation*}
\begin{equation*}\label{}
[x_{i,r+1}^{\pm}, x_{j,s}^{\pm}]-
[x_{i,r}^{\pm}, x_{j,s+1}^{\pm}]=\pm\frac12
d_ia_{ij}\left(x_{i,r}^{\pm}x_{j,s}^{\pm}
+x_{j,s}^{\pm}x_{i,r}^{\pm}\right),
\end{equation*}
\begin{equation*}\label{}
\sum_\pi
[x_{i,r_{\pi\left(1\right)}}^{\pm},
[x_{i,r_{\pi\left(2\right)}}^{\pm}, \ldots,
[x_{i,r_{\pi\left(m\right)}}^{\pm},
x_{j,s}^{\pm}]\cdots]]=0, i\neq j,
\end{equation*}
for all sequences of non-negative integers $r_1,\ldots,r_m$, where
$m=1-a_{ij}$ and the sum is over all permutations $\pi$ of $\{1,\dots,m\}$.
\end{definition}

Let the degree of $x_{i,d}^{\pm}$ and $h_{i,d}$ be $d$; then we obtain a filtration on $\yg$ by taking $\yg_r$ to be the linear span of all monomials of degree at most $r$ in the generators $x_{i,s}^{\pm}$ and $h_{i,s}$.

\begin{proposition}[Proposition 12.1.6, \cite{ChPr2}]
The associated graded algebra $\operatorname{gr}\yg$ is isomorphic to the universal enveloping algebra of the current algebra $\g\otimes_\C\C[t]$.
\end{proposition}
There is a version of the Poincare-Birkhoff-Witt theorem for Yangians.

\begin{proposition}[\cite{Le}] Let $\Delta^{+}$ be the set of positive roots of $\g$. Let $x_{\alpha,k}^{\pm}$ be the elements constructed in \cite{Le}.
Fix a total ordering on the set $$\sum=\{x_{\alpha,k}^{\pm}\ |\ \alpha\in \Delta^{+}, r\in\Z_{\geq 0}\}\ \bigcup\  \{h_{i,r}\ |\ i\in I, r\in\Z_{\geq 0}\}.$$ Then the set of ordered monomials in the elements of $\sum$ is a vector space basis of $\yg$.
\end{proposition}

It has been established that $\yg$ has a Hopf algebra structure. The coproduct of $\yg$, however, is not given explicitly in terms of the generators $x_{i,r}^{\pm}$ and $h_{i,r}$.
To describe the Hopf algebra structure of $\yg$, let $Y^{\pm}$ and $H$ be the subalgebras of $Y\left(\mathfrak{g}\right)$ generated by the $x_{\alpha, k}^{\pm}$, for all positive roots $\alpha$, and $h_{i, k}$ for all $i\in I$, respectively.
For a fixed $i\in I$, let $Y_{i}^{\pm}$ and $H_i$ be the subalgebras of $\yg$ generated by $x_{i,k}^{\pm}$ and $h_{i,k}$, respectively, and set $M_{i}^{\pm}=\sum_{k} x_{i, k}^{\pm}Y_{i}^{\pm}$. Set $N^{\pm}=\sum_{i,k} x_{i, k}^{\pm}Y^{\pm}$ and define $N^{\pm}_{i}$ to be the subalgebras of $\yg$ generated by all monomials in $x_{\alpha,k}^{\pm}$ with at least one factor with $\alpha\neq\alpha_i$.  V. Chari and A. Pressley proved the following proposition.
\begin{proposition}[Proposition 2.8,\cite{ChPr4}]\label{Delta}\ \\
Denote $\yg$ by $Y$. Modulo $\overline{Y}\equiv \left(N^{-}Y\otimes YN^{+}\right)\bigcap \left(YN^{-}\otimes N^{+}Y\right)$, we have
\begin{enumerate}
  \item $\Delta_{\yg}\left(x_{i,k}^{+}\right)\equiv x_{i,k}^{+}\otimes 1+1\otimes x_{i,k}^{+}+\sum\limits_{j=1}^{k}h_{i,j-1}\otimes x_{i,k-j}^{+}$,
  \item $\Delta_{\yg}\left(x_{i,k}^{-}\right)\equiv x_{i,k}^{-}\otimes 1+1\otimes x_{i,k}^{-}+\sum\limits_{j=1}^{k}x_{i,k-j}^{-}\otimes h_{i,j-1}$,
  \item $\Delta_{\yg}\left(h_{i,k}\right)\equiv h_{i,k}\otimes 1+1\otimes h_{i,k}+\sum\limits_{j=1}^{k}h_{i,j-1}\otimes h_{i,k-j}$.
\end{enumerate}
\end{proposition}
We need to modify some of the results of Proposition \ref{Delta}. We denote the generators of $\ysl$ by $x_{k}^{\pm}$ and $h_{k}$. In $\yg$, for a fixed $i\in I$, $\{x_{i,r}^{\pm}, h_{i,r}|r\geq 0\}$ generates a subalgebra $Y_i$ of $\yg$ and the assignments $\frac{\sqrt{d_i}}{{d_i}^{k+1}}x_{i,k}^{\pm}\rightarrow x_{k}^{\pm}$ and $\frac{1}{{d_i}^{k+1}}h_{i,k}\rightarrow h_{k}$ define an isomorphism from $Y_i$ to $\ysl$. Let $\overline{h}_{1}=h_{1}-\frac{1}{2}h_{0}^2$ and $\overline{h}_{i,1}=h_{i,1}-\frac{1}{2}h_{i,0}^2$. It is proved in \cite{ChPr5} that $\Delta_{\ysl}\left(\bar{h}_{1}\right)=\bar{h}_{1}\otimes 1+1\otimes \bar{h}_{1}-2x^{-}_{0}\otimes x^{+}_{0}.$ The isomorphism $Y_i\cong \ysl$ becomes a Hopf algebra isomorphism if we set $\Delta_{Y_i}\left(\bar{h}_{i,1}\right)=\bar{h}_{i,1}\otimes 1+1\otimes \bar{h}_{i,1}-\left(\alpha_i,\alpha_i\right)x^{-}_{i,0}\otimes x^{+}_{i,0}.$  We denote by $\Delta_{Y_i}$ and $\Delta_{\yg}$ the coproducts of $Y_i$ and $\yg$, respectively. There is a difference between the coproduct of $Y_i$ and the coproduct of $\yg$ restricted to $Y_i$. For instance,
$$\Delta_{Y_i}\left(\bar{h}_{i,1}\right)=\bar{h}_{i,1}\otimes 1+1\otimes \bar{h}_{i,1}-\left(\alpha_i,\alpha_i\right)x^{-}_{i,0}\otimes x^{+}_{i,0};$$
$$\Delta_{\yg}\left(\bar{h}_{i,1}\right)=\bar{h}_{i,1}\otimes 1+1\otimes \bar{h}_{i,1}-\left(\alpha_i,\alpha_i\right)x^{-}_{i,0}\otimes x^{+}_{i,0}-\sum\limits_{\alpha\succ 0 , \alpha\neq \alpha_i}\left(\alpha,\alpha_i\right)x^{-}_{\alpha,0}\otimes x^{+}_{\alpha,0}.$$

The next proposition follows from the proof of Proposition 2.8 \cite{ChPr4}.

\begin{proposition}
 $\Delta_{\yg}\left(x^{\pm}_{i,k}\right)- \Delta_{Y_i}\left(x^{\pm}_{i,k}\right)\in H_iN^{-}_i\otimes H_iN^{+}_i$.
\end{proposition}
Similar arguments show that
\begin{proposition}\label{c1dgd2d}\
\begin{enumerate}
  \item $\Delta_{\yg}\left(x_{i, k_s}^{-}\ldots x_{i, k_1}^{-}\right)-\Delta_{Y_i}\left(x_{i, k_s}^{-}\ldots x_{i, k_1}^{-}\right)\in HN^{-}_{i}\otimes M_i^{-}HN^{+}_{i}$.
  \item $\Delta_{\yg}\left(x_{i, k_s}^{+}\ldots x_{i, k_1}^{+}\right)-\Delta_{Y_i}\left(x_{i, k_s}^{+}\ldots x_{i, k_1}^{+}\right)\in M_i^{+}HN^{-}_{i}\otimes HN^{+}_{i}$.
\end{enumerate}

\end{proposition}
Before  moving on to discuss the representation theory of $\yg$, we present the following useful proposition.
\begin{proposition}[Proposition 2.6, \cite{ChPr4}]\label{s1.2.1p26} For any $a\in \C$,
the assignment $$\tau_{a}(h_{i,k})=\sum_{r=0}^{k}\binom{k}{r}a^{k-r}h_{i,r},\qquad \tau_{a}(x_{i,k}^{\pm})=\sum_{r=0}^{k}\binom{k}{r}a^{k-r}x_{i,r}^{\pm}$$ extends to a Hopf algebra automorphism of $\yg$.
\end{proposition}

We are now in the position to introduce basic results on the finite-dimensional representation theory of $\yg$.
\begin{definition}
A representation $V$ of the Yangian $Y\left(\mathfrak{g}\right)$ is said to be highest weight if it is generated by a vector $v^{+}$ such that $x_{i,k}^{+}v^{+}= 0$ and $h_{i,k}v^{+} =\mu_{i,k} v^{+}$, where $\mu_{i,k}$ are complex numbers.
\end{definition}

The defining relations of a Yangian allow one to define highest weight representations of $\yg$. Set
$\mu=\Big(\mu_1\left(u\right),\mu_2\left(u\right),\ldots,\mu_l(u)\Big)$,
where $\mu_i\left(u\right)=1+\sum\limits_{k\in \Z_{\geq 0}}\mu_{i,k}u^{-k-1}$ is a formal series in $u^{-1}$ for $i\in I$.
The Verma module $M\left(\mu\right)$ of $\yg$ is defined to be the quotient of $\yg$ by the left ideal generated by $N^{+}$
and the elements $h_{i,k}-\mu_{i,k}1$. The image $1_\mu$ of the element $1\in \yg$ in the quotient is a highest weight
vector of $M\left(\mu\right)$. The Verma module $M(\mu)$ is a universal highest weight representation in the sense that if
$V\left(\mu\right)$ is another highest weight representation with a highest weight vector $v$, then the mapping
$1_{\mu}\mapsto v$ defines a surjective $\yg$-module homomorphism $M\left(\mu\right)\rightarrow V\left(\mu\right)$.
It is known that $\yg_0$, the degree 0 part of the filtration on $\yg$, is isomorphic to the universal enveloping algebra
$U(\g)$, hence $V(\mu)$ can be viewed as a $U(\g)$-module. The weight subspace $V_{\mu^{\left(0\right)}}$ with
$\mu^{\left(0\right)}=\left(\mu_{1,0},\dots,\mu_{l,0}\right)$ is one-dimensional and spanned
by the highest weight vector of $V(\mu)$. All other nonzero weight subspaces
correspond to weights $\gamma$ of the form
$\gamma=\mu^{\left(0\right)}-k_1\alpha_1-\ldots-k_l\alpha_l,$
where all $k_i$ are nonnegative integers and not all of them are zero.
A standard argument shows that $M\left(\mu\right)$ has a unique irreducible quotient $L\left(\mu\right)$.
In \cite{Dr2}, V. Drinfeld gave the following classification of the finite-dimensional irreducible representations of $\yg$.
\begin{theorem}\label{cfdihwr}\
\begin{enumerate}
  \item[(a)] Every irreducible finite-dimensional representation of $\yg$ is highest weight.
  \item [(b)] The irreducible representation $L\left(\mu\right)$ is finite-dimensional if and only if there exists an $l$-tuple of polynomials $\pi=\left(\pi_1(u),\ldots, \pi_l(u)\right)$ such that $$\frac{\pi_i\left(u+d_i\right)}{\pi_i\left(u\right)}=1+\sum_{k=0}^{\infty} \mu_{i,k}u^{-k-1},$$ in the sense that the right-hand side is the Laurent expansion of the left-hand side about $u=\infty$. The polynomials $\pi_i\left(u\right)$ are called Drinfeld polynomials.
\end{enumerate}
\end{theorem}

In \cite{ChPr4}, V. Chari and A. Pressley proved that $L(\mu)$ is a subquotient of a tensor product of fundamental representations, see Theorem \ref{efdiriasqotp}. To show this, we need the following proposition.

\begin{proposition}[Proposition 2.15, \cite{ChPr4}]\label{vtv'hwv}
Let both $V$ and $\widetilde{V}$ be irreducible finite-dimensional representations of $Y\left(\mathfrak{g}\right)$ with associated l-tuples of polynomials $\pi$ and $\widetilde{\pi}$, respectively. Let $v^{+}\in V$, $\widetilde{v}^{+}\in \widetilde{V}$ be their highest weight vectors. Then $v^{+}\otimes \widetilde{v}^{+}$ is a highest weight vector in $V\otimes \widetilde{V}$ and its associated polynomials are $\pi_i\widetilde{\pi}_i$.
\end{proposition}
\begin{theorem}[Theorem2.16, \cite{ChPr4}]\label{efdiriasqotp}
 Every finite-dimensional irreducible representation $L(\mu)$ of $\yg$ is a subquotient of a tensor product $W=V_1\otimes\ldots\otimes V_n$ of the fundamental representations. In fact, $V$ is a quotient of the cyclic sub-representation of $W$ generated by the tensor product of the highest weight vectors in the $V_i$.
\end{theorem}
In \cite{ChPr3} and \cite{ChPr5}, it was proved that every finite-dimensional irreducible representation of $\ysl$ is a tensor product of evaluation representations  which are irreducible under $\mathfrak{sl}_2$. However, this property is not true in general and in most cases the structure of finite-dimensional irreducible $\yg$-modules remains unknown.

Let $V$ be a finite-dimensional irreducible representation of $\yg$ with $l$-tuple of Drinfeld polynomials $\pi$. Define the following associated $\yg$-representations:

\begin{enumerate}
  \item Pulling back $V$ through $\tau_{a}$ as defined in Proposition \ref{s1.2.1p26}, we denote the representation by $V(a)$. $V(a)$ has Drinfeld polynomials $\big(\pi_1(u-a),\ldots, \pi_l(u-a)\big)$.
  \item The left dual $^t V$ of $V$ and right dual $V^t$ of $V$ are the representations of $\yg$ on the dual vector space of $V$ defined as follows:
\begin{center}
$\left(y\cdot f\right)\left(v\right)=f\left(S\left(y\right)\cdot v\right)$, \qquad $y\in \yg$, $f\in\ ^tV$, $v\in V$;\\
$\left(y\cdot f\right)\left(v\right)=f\left(S^{-1}\left(y\right)\cdot v\right)$, \qquad $y\in \yg$, $f\in V^t$, $v\in V$,
\end{center}
where $S$ is the antipode of $\yg$.
\end{enumerate}
It is well known that the dual of an irreducible representation is irreducible.
The dual spaces $\ ^tV$ and $V^t$ are vital in determining the irreducibility of $V$ as indicated in the proposition below.
\begin{proposition}[Proposition 3.8, \cite{ChPr8}]\label{VoWWoVhi}
Let $V$ be a finite-dimensional representation of $\yg$. $V$ is irreducible if and only if $V$ and $\ ^tV$ \big(respectively, $V$ and $V^t$\big) are both highest weight $\yg$-modules.
\end{proposition}
Note that every finite-dimensional highest-weight representation is also a lowest weight representation and vice-versa. The following lemma is an analogue of Lemma 4.2 in \cite{Ch3}. The proof has been omitted as it is similar to the proof of that Lemma.
\begin{lemma}\label{v-w+gvtw}
Let $V$ and $W$ be two finite-dimensional highest weight representations of $\yg$ with lowest and highest weight vectors $v^{-}$ and $w^{+}$, respectively. Then $v^{-}\otimes w^{+}$generates $V\otimes W$.
\end{lemma}
\section{Local Weyl modules of $\yg$ have finite dimension}

In this section, we give the definition of the local Weyl module $W(\pi)$ of $\yg$ and show that $W(\pi)$ is finite-dimensional.
\begin{definition}\label{TFPD31}
Let $\pi_i\left(u\right)=\prod\limits_{j=1}^{m_i}\left(u-a_{i,j}\right)$ and $\pi=\big(\pi_1(u),\pi_2(u),\ldots, \pi_l(u)\big)$. The local Weyl module $W\left(\mathbf{\pi}\right)$ for $\yg$ is defined as the module generated by a highest weight vector $w_\pi$ that satisfies the following relations: \begin{equation}
x_{i,k}^{+}w_\pi=0, \qquad \left(x_{i,0}^{-}\right)^{m_i+1}w_\pi=0,\qquad \left(h_{i}\left(u\right)-\frac{\pi_i\left(u+d_i\right)}{\pi_i\left(u\right)}\right)w_\pi=0.\end{equation}
\end{definition}

\begin{remark}\label{mi=lambdahi}\
\begin{enumerate}
  \item The local Weyl module $W(\pi)$ of $\yg$ can be considered as the quotient of the Verma module $M(\mu)$ by the submodule generated by $\left(x_{i,0}^{-}\right)^{m_i+1}w_\pi$, where $\mu_i(u)=\frac{\pi_i(u+d_i)}{\pi_i(u)}$.
  \item It follows from Remark [2], part C of section 12.1 of \cite{ChPr2} that for any finite-dimensional highest weight representation $V(\pi)$ of $\yg$, $\left(x_{i,0}^{-}\right)^{m_i+1}w_\pi=0$. Therefore $V(\pi)$ is a quotient of $W(\pi)$.
\end{enumerate}
\end{remark}
Before showing that $W(\pi)$ is finite-dimensional, we need some results regarding local Weyl modules of the current algebra $\g[t]$.
\begin{definition}\label{c1lwmotca}
Let $\lambda=\sum\limits_{i\in I} m_i \omega_i$ be a dominant integral weight of $\g$. Denote by $W\left(\lambda\right)$ the $\g[t]$-module generated by an element $v_{\lambda}$ which satisfies the relations:
$$
\mathfrak{n}^{+}\otimes \C[t] v_{\lambda} = 0, \;h \otimes t\C[t] v_{\lambda} = 0, \; hv_{\lambda} = \lambda\left(h\right)v_{\lambda},\ \text{and}\ \  \left(x_{\alpha_i}^{-}\otimes 1\right)^{m_i +1}v_{\lambda}= 0,
$$
for all $h \in \mathfrak{h}$ and all simple roots $\alpha_i$. This module is called
the local Weyl module associated to $\lambda\in P^+$.
\end{definition}

\begin{theorem}[Theorem 1.2.2, \cite{ChPr1}]\label{mpocac1}
Let $\lambda$ be a dominant integral weight of $\g$. The local Weyl module $W\left(\lambda\right)$ is finite-dimensional. Moreover, any finite-dimensional $\g[t]$-module $V$ generated by an element $v\in V$ satisfying the relations:
$$
\mathfrak{n}^{+}\otimes \C[t] v = 0, \quad h \otimes t\C[t] v = 0, \quad hv = \lambda\left(h\right)v
$$ is a quotient of $W\left(\lambda\right)$.
\end{theorem}
The dimension of the local Weyl module $W(\lambda)$ of $\g[t]$ has been determined.
\begin{proposition}[Corollary A, \cite{Na}]\label{PS3TADOLWM}
$$\operatorname{Dim}\ W\left(\lambda\right)=\prod_{i\in I} \Big(\operatorname{Dim}\big(W(\omega_i)\big)\Big)^{m_i}.$$
\end{proposition}
Combining the isomorphism (9.1) in \cite{Na} with the first part of the main theorem of Section 2.2 in \cite{ChMo}, we obtain the following important corollary.
\begin{corollary}\label{dkrvawocsp} When $\g$ is a classical simple Lie algebra, $W\left(\omega_i\right)\cong_{\g} V_a\left(\omega_i\right)$ for any $i\in I$ and $a\in \C$. In particular,
$\operatorname{Dim}\big(V_a(\omega_i)\big)=\operatorname{Dim}\big(W(\omega_i)\big)$.
\end{corollary}
We are now in a position to prove the main objective of this section.
\begin{theorem}\label{wmifd}
The local Weyl module $W(\pi)$ is finite-dimensional.
\end{theorem}
\begin{proof} Let $\lambda=\sum\limits_{i\in I}m_i\omega_i$, where $m_i$ is the degree of the polynomial $\pi_i\left(u\right)$.
We divide the proof into steps.

Step 1: The $\N$-filtration on $\yg$ induces a filtration on $W(\pi)$ as follows: let $W(\pi)_s$ be the subspace of $W(\pi)$ spanned by elements of the form $yw_{\pi}$, where $y\in \yg_{s}$.
Denote by gr$\big(W(\pi)\big)=\bigoplus\limits_{r=0}^{\infty} W(\pi)_{r}/W(\pi)_{r-1}$ the associated graded module, where $W(\pi)_{-1}=0$. Let $\overline{w_{\pi}}$ be the image of $w_{\pi}$ in gr$\big(W(\pi)\big)$.



Step 2: $\text{gr}\big(W(\pi)\big)=\text{gr}\big(\yg\big)
\overline{w_{\pi}}=U\left(\g[t]\right)\overline{w_{\pi}}$.

Proof:  Let $\overline{v}\in \text{gr}\big(W(\pi)\big)$. Without loss of generality, we may assume that $v\in W(\pi)_r\backslash W(\pi)_{r-1}.$ Thus $v=yw_{\pi}$ for some $y\in \yg_r\backslash \yg_{r-1}$ by Step 1 and hence $\bar{v}=\bar{y}\ \overline{w_{\pi}},$ where $\bar{y}$ is the image of $y$ in $\yg_r/\yg_{r-1}.$  So $\text{gr}\big(W(\pi)\big)=\text{gr}\big(\yg\big)\overline{w_{\pi}}.$


Step 3: $\text{gr}\big(W(\pi)\big)$ is a highest weight representation of $\text{gr}\big(\yg\big)$.

Proof: Note that $\left(x_{\alpha}^{+}\otimes t^r\right)\overline{w_{\pi}}=\overline{x_{\alpha,r}^{+}w_{\pi}}= \overline{0}$ and $h_i\overline{w_{\pi}}=\overline{\frac{1}{d_i}h_{i,0}w_{\pi}}=m_i\overline{w_{\pi}}$. Since $w_\pi\in W(\pi)_0\subseteq W(\pi)_{r}$ for all $r\geq 1$, $\left(h_i\otimes t^r\right)\overline{w_{\pi}}=\overline{h_{i,r}w_{\pi}}=
\overline{c.w_{\pi}}=\overline{0}$, where $c\in\C$.  Moreover, $(x_{i}^-)^{m_i+1}\overline{w_{\pi}}=\overline{(x_{i,0}^-)^{m_i+1}w_{\pi}}=0$. Therefore Step 3 holds.

Step 4: $W(\pi)$ is finite-dimensional.

Proof: There is a surjective homomorphism $\varphi: W\left(\lambda\right)\rightarrow \text{gr}\big(W(\pi)\big)$ by Theorem \ref{mpocac1}. It follows from Proposition \ref{PS3TADOLWM} that $W(\lambda)$ is finite-dimensional and hence $\text{gr}\big(W(\pi)\big)$ is as well since it is a quotient of $W(\lambda)$. It follows that  $W(\pi)$ is finite-dimensional and $\operatorname{Dim}\big(W(\pi)\big)= \operatorname{Dim}\Big(\text{gr}\big(W(\pi)\big)\Big)$.
\end{proof}

Since $\text{gr}\big(W(\pi)\big)$ is a quotient of the Weyl module $W\left(\lambda\right)$ of $\g[t]$, we obtain an upper bound of the dimension of the local Weyl module.
\begin{theorem}\label{ubodowm} Using the notation as in the theorem above,
$$\text{Dim}\big(W(\pi)\big)\leq\text{Dim}\big(W(\lambda)\big).$$
\end{theorem}

\section{Local Weyl module $W(\pi)$ of $\ysl$}
In this section, we describe the local Weyl module $W(\pi)$ of $\ysl$. This module plays a key role in the  characterization of local Weyl modules of $\yg$.

Denote by $W_m(a)$ the finite-dimensional irreducible representation of $\ysl$ associated to the Drinfeld polynomial $\big(u-a\big)\big(u-(a_1+1)\big)\ldots \big(u-(a+m-1)\big)$.

\begin{proposition}[Proposition 3.5,\cite{ChPr3}]\label{wra}
For any $m\geq 1$, the module $W_m(a)$ has a basis
$\{w_0,w_1,\ldots, w_m\}$ on which the action of $\ysl$ is given by
\begin{center}
$ x_k^+w_s=(s+a)^k(s+1)w_{s+1},\ \ x_k^-w_s=(s+a-1)^k(m-s+1)w_{s-1},$\\

$h_kw_s=\big((s+a-1)^ks(m-s+1)-(s+a)^k(s+1)(m-s)\big)w_s$.
\end{center}
\end{proposition}
The cases when $m=1$ and $m=2$ are important in the characterization of local Weyl modules of $\yg$. Explicit identities are listed in the corollary below.
\begin{corollary}\label{slw1a}\
In $W_1(a)$,
\begin{equation}\label{Cor3.3}
  h_{k}w_1=a^kw_1,\quad x_{k}^{-}w_1=a^kx_{0}^{-}w_1,\quad h_{k}x_{i,0}^{-}w_1=-a^kx_{0}^{-}w_1.
\end{equation}

In $W_2(a)=\ysl(w_2)=span\{w_0, w_1, w_2\}$,
\begin{enumerate}
    \item $(x_{1}^{-}x_{0}^{-}+x_{0}^{-}x_{1}^{-})(w_2)=(2a+1)(x_{0}^{-})^2(w_2)$.
    \item $(x_{2}^{-}x_{0}^{-}+x_{0}^{-}x_{2}^{-})(w_2)=(2a^2+2a+1)(x_{0}^{-})^2(w_2)$.
\end{enumerate}
\end{corollary}
It follows from Proposition \ref{Delta} and Corollary \ref{slw1a} that:
\begin{corollary}\label{w1bw1a}
Let $v_1$ and $w_1$ be the highest weight vectors in $W_1(b)$ and $W_1(a)$ respectively.
In $W_1(b)\otimes W_1(a)$,
\begin{enumerate}
  \item $(x_{1}^{-}x_{0}^{-}+x_{0}^{-}x_{1}^{-})(v_1\otimes w_1)=(a+b)(x_{0}^{-})^2(v_1\otimes w_1)$.
  \item 
      $(x_0^{-}x_{2}^{-}+x_{2}^{-}x_0^{-})(v_1\otimes w_1)=(a^2+b^2)(x_{0}^{-})^2(v_1\otimes w_1)$.

\end{enumerate}
\end{corollary}

The next proposition is fundamental for the characterization of finite-dimensional irreducible representations of $\ysl$.

\begin{proposition}[Proposition 3.7, \cite{ChPr3}]\label{ctpihw}
Let $a_1,a_2,\ldots,a_m\in\C$, $m\geq 1$. Then if $a_j-a_i\neq 1$ when $i<j$,  $W_1(a_1)\otimes\ldots\otimes W_1(a_m)$ is a highest weight $\ysl$-module.
\end{proposition}

From this proposition, we obtain:
\begin{corollary}\label{tpihw}
Let $a_1,a_2,\ldots,a_m\in\C$, $m\geq 1$, and $\operatorname{Re}(a_1)\geq \ldots \geq \operatorname{Re}(a_m)$, where $\operatorname{Re}(a_i)$ is the real part of $a_i$. Then $W_1(a_1)\otimes\ldots\otimes W_1(a_m)$ is a highest weight $\ysl$-module.
\end{corollary}

\begin{lemma}[Corollary 3.8.\cite{ChPr3}]\label{beautiful}
Let $P_1, P_2,\ldots, P_m$ be polynomials, and let $V(P_i)$ be the irreducible highest weight representation whose Drinfeld polynomial is $P_i$. Assume that if $a_i$ is a root of $P_i$ and $a_j$ a root of $P_j$, where $i<j$, then $a_j-a_i\neq 1$. Then $V(P_1)\otimes V(P_2)\otimes \ldots\otimes V(P_m)$ is a highest weight $\ysl$-module.
\end{lemma}

\begin{theorem}\label{wmfsl2} Let $\pi(u)=(u-a_1)(u-a_2)\ldots(u-a_m)$ be a decomposition of $\pi$ over $\C$ such that $\operatorname{Re}(a_1)\geq \ldots \geq \operatorname{Re}(a_m)$.
Then $W(\pi)$ is isomorphic to the ordered tensor product $W_1(a_1)\otimes\ldots\otimes W_1(a_m)$.
\end{theorem}
\begin{proof}
It follows from Corollary \ref{tpihw} that $W_1(a_1)\otimes\ldots\otimes W_1(a_m)$ is a highest weight module of $\ysl$ and from Proposition \ref{vtv'hwv} we obtain that the associated polynomial is $\pi(u)$. As a consequence of Remark \ref{mi=lambdahi}, we have $\text{Dim}\big(W(\pi)\big)\geq 2^m$. By the main theorem of \cite{ChLo} and Theorem \ref{ubodowm},
$\text{Dim}\big(W(\pi)\big)\leq 2^m$.
Thus $\text{Dim}\big(W(\pi)\big)=2^m$. Therefore $W(\pi)=W_1(a_1)\otimes\ldots\otimes W_1(a_m).$
\end{proof}
\begin{remark}
The order of the tensor product decomposition of $W(\pi)$ in the previous Theorem is not necessarily unique. We use the condition $\operatorname{Re}(a_1)\geq \ldots \geq \operatorname{Re}(a_m)$ because it is convenient.
\end{remark}
\section{Local Weyl module $W(\pi)$ of $\yg$}
Let $v_{1}^{+}$ and $v_{1}^{-}$ be the highest and lowest weight vectors in $V_a(\omega_i)$, respectively. There exists a ``path" from  $v_{1}^{+}$ to $v_{1}^{-}$. For instance, when $\g=\nyn$,
$$v^{-}_1=\left(x_{l+1-i,0}^{-}\ldots x_{2,0}^{-}x_{1,0}^{-}\right)\left(x_{l-i+2,0}^{-}\ldots x_{3,0}^{-}x_{2,0}^{-}\right)\ldots \left(x_{l,0}^{-}\ldots x_{i+1,0}^{-}x_{i,0}^{-}\right)v_1^{+}.$$
When $\g$ is a simple Lie algebra of other than type $A$, we use the techniques in \cite{Ch3} to find a path. Suppose that one reduced expression of the longest element of the Weyl group of $\g$ is $w_0=s_{r_1}s_{r_2}\ldots s_{r_p}$, where $s_{r_j}$ is a simple reflection. Suppose $s_{r_{j+1}}s_{r_{j+2}}\ldots s_{r_p}(\omega_i)=m_j\omega_{r_j}+\sum\limits_{n\neq r_j}c_n\omega_n$. Then
\begin{center}
$v_{1}^{-}=\big(x_{r_1,0}^{-}\big)^{m_{1}}\big(x_{r_2,0}^{-}\big)^{m_{2}}\ldots \big(x_{r_p,0}^{-}\big)^{m_{p}}v_{1}^{+}.$
\end{center}
Define $\sigma_j:=s_{r_{j+1}}s_{r_{j+2}}\ldots s_{r_p}$ and $v_{\sigma_j(\omega_i)}:=\big(x_{r_{j+1},0}^{-}\big)^{m_{j+1}}\big(x_{r_{j+2},0}^{-}\big)^{m_{j+2}}\ldots \big(x_{r_p,0}^{-}\big)^{m_{p}}v_{1}^{+}$.
We claim that
\begin{lemma}\label{S5L1yiahwr} Suppose $m_j\neq 0$.
$Y_{r_j}\left(v_{\sigma_j(\omega_i)}\right)$ is a highest weight representation of $Y_{r_j}$.
\end{lemma}
\begin{proof}
Suppose that $x_{r_j,0}^{+}v_{\sigma_j(\omega_i)}\neq 0$,
then $\sigma_j(\omega_{i})+\alpha_{r_j}$ is a weight of $V_{a_1}(\omega_i)$.    Thus $\omega_i+\sigma_{j}^{-1}(\alpha_{r_j})$ is a weight which does not precede $\omega_i$, which implies that $\sigma_{j}^{-1}(\alpha_{r_j})$ is a negative root.
However, $\ell(s_{r_j}\sigma_j)=\ell(\sigma_j)+1$,  hence $\sigma_{j}^{-1}(\alpha_{r_j})$ is a positive root, so we have a contradiction.
Thus $x_{r_j,0}^{+}v_{\sigma_j(\omega_i)}=0$.
Since the weight space of weight $\sigma_j(\omega_i)$ is 1-dimensional and $H$ is a commutative subalgebra of $\yg$, $h_{r_j,s}v_{\sigma_j(\omega_i)}=c_{j,s}v_{\sigma_j(\omega_i)}$, where $c_{j,s}\in \C$.
Thus $Y_{r_j}\left(v_{\sigma_j(\omega_i)}\right)$ is a highest weight representation of $Y_{r_j}$.
\end{proof}
As $V_{a}(\omega_i)$ is finite-dimensional, so is $Y_{r_j}(v_{\sigma_j(\omega_i)})$. Let the associated polynomial of $Y_{r_j}\left(v_{\sigma_j(\omega_i)}\right)$ be $Q_{i,j}(u)$ and set $Q_{i,j}:=Q_{i,j}(u)$.
The highest weight representation $Y_{r_j}\left(v_{\sigma_j(\omega_i)}\right)$ can be partially understood via its associated polynomial $Q_{i,j}$: it is a quotient of the local Weyl module $W(Q_{i,j})$ of $\ysl$.

\begin{theorem}\label{S5T1Liahw}
The ordered tensor product $L=V_{a_1}(\omega_{b_1})\otimes V_{a_2}(\omega_{b_2})\otimes\ldots\otimes V_{a_k}(\omega_{b_k})$ is a highest weight representation if for all $1\leq j\leq p$ and $1\leq m<n\leq k$, when $b_n=r_j$, the difference of the number $\frac{a_n}{d_{r_j}}$ and any root of the associated polynomial of $Y_{r_j}\left(v_{\sigma_j(\omega_{b_m})}\right)$ does not equal 1.
\end{theorem}
\begin{proof}
Let $v_m^{+}$ be the highest weight vector of the fundamental module $V_{a_m}\left(\omega_{b_m}\right)$ and let $v_1^{-}$ be the lowest weight vector of $V_{a_1}\left(\omega_{b_1}\right)$. We prove this theorem by induction on $k$.
For $k=1$, $L=V_{a_{1}}\left(\omega_{b_{1}}\right)$ is irreducible and hence it is a highest weight representation. We assume that the claim is true for all positive integers less than or equal to $k-1$ ($k\geq 2$). By the induction hypothesis, $V_{a_2}(\omega_{b_2})\otimes V_{a_3}(\omega_{b_3})\otimes\ldots\otimes V_{a_k}(\omega_{b_k})$ is a highest weight representation of $\yg$ and its highest weight vector is $v^{+}=v_2^{+}\otimes \ldots\otimes v_k^{+}$. To show that $L$ is a highest weight representation, it suffices to show $v^{-}_{1}\otimes v^{+}\in \yg\left(v^{+}_1\otimes v^{+}\right)$ by Lemma \ref{v-w+gvtw}.

We first show that $$Y_{r_j}\left(v_{\sigma_j\left(\omega_{b_1}\right)}\right)\otimes Y_{r_j}\left(v_2^{+}\right)\otimes\ldots\otimes Y_{r_j}\left(v_k^{+}\right)= Y_{r_j}\left(v_{\sigma_j\left(\omega_{b_1}\right)}\otimes v_2^{+}\otimes\ldots \otimes v_k^{+}\right).$$
Recall that we normalized the generators of $Y_{r_j}$ to satisfy the defining relations of $\ysl$. For $2\leq n\leq k$, $Y_{r_j}(v_n^{+})$ is either trivial, if $r_j\neq b_n$, or isomorphic to $W_1(\frac{a_n}{d_{r_j}})$, if $r_j=b_n$. When  $r_j=b_n$, if the difference of the number $\frac{a_n}{d_{r_j}}$ and any root of $Q_{b_1,j}$ does not equal 1, $W(Q_{b_1,j})\otimes Y_{r_j}\left(v_2^{+}\right)\otimes\ldots\otimes Y_{r_j}\left(v_k^{+}\right)$ is a highest weight representation of $\ysl$ by Lemma \ref{beautiful}.
Hence its quotient $Y_{r_j}\left(v_{\sigma_j\left(\omega_{b_1}\right)}\right)\otimes Y_{r_j}\left(v_2^{+}\right)\otimes\ldots\otimes Y_{r_j}\left(v_k^{+}\right)$ is a highest weight module of $Y_{r_j}$ with highest weight vector $v_{\sigma_j\left(\omega_{b_1}\right)}\otimes v_2^{+}\otimes\ldots \otimes v_k^{+}$. Thus $$Y_{r_j}\left(v_{\sigma_j\left(\omega_{b_1}\right)}\otimes v_2^{+}\otimes\ldots \otimes v_k^{+}\right)\supseteq Y_{r_j}\left(v_{\sigma_j\left(\omega_{b_1}\right)}\right)\otimes Y_{r_j}\left(v_2^{+}\right)\otimes\ldots\otimes Y_{r_j}\left(v_k^{+}\right).$$
By the coproduct of Yangians and Proposition \ref{c1dgd2d}, it is obvious that $$Y_{r_j}\left(v_{\sigma_j\left(\omega_{b_1}\right)}\otimes v_2^{+}\otimes\ldots \otimes v_k^{+}\right)\subseteq Y_{r_j}\left(v_{\sigma_j\left(\omega_{b_1}\right)}\right)\otimes Y_{r_j}\left(v_2^{+}\right)\otimes\ldots\otimes Y_{r_j}\left(v_k^{+}\right).$$ Therefore the claim is true.

Since $v_{\sigma_{j-1}(\omega_{b_1})}\otimes v^{+}=(x_{r_j,0}^{-})^{m_j}v_{\sigma_j\left(\omega_{b_1}\right)}\otimes v^{+}\in Y_{r_j}\left(v_{\sigma_{j}(\omega_{b_1})}\right)\otimes Y_{r_j}\left(v^{+}\right)$, we have $v_{\sigma_{j-1}(\omega_{b_1})}\otimes v^{+}\in Y_{r_j}\left(v_{\sigma_j\left(\omega_{b_1}\right)}\otimes v^{+}\right)$.
By downward induction on the subscript $j$ of $v_{\sigma_{j}(\omega_{b_1})}$,
$v^{-}_{1}\otimes v^{+}\in \yg\left(v^{+}_1\otimes v^{+}\right)$. By Lemma \ref{v-w+gvtw}, $L=\yg \left(v_1^{+}\otimes v^{+}\right)$.
\end{proof}

To find an explicit cyclicity condition for the ordered tensor product $L$, we need to find the roots of the polynomial $Q_{b_m,j}$. From now on, we suppose that $\g$ is a classical simple Lie algebra over $\C$. In contrast to the method used in Proposition 6.3 in \cite{Ch3}, we compute the associated polynomial $Q_{b_m,j}$ by using some of the defining relations of $\yg$. We only provide the details for the case when $\g=\nysp$ where all the techniques needed to prove the other cases are used (the detailed computations in the other cases can be found in \cite{Ta}). Before carrying out some long computations, we would like to indicate to our reader that a concrete  cyclicity condition for the ordered tensor product $L$ is given in Theorem \ref{mtoyspl}.
\begin{remark}
In this paper, we fix the reduced expression for the longest element $w_0$ given in table 1 in \cite{BeKaOhPa}. When $\g=\nyo$,
\begin{align*}
  w_0&=s_{l}s_{l-1}\left(s_{l-2}s_{l}s_{l-1}s_{l-2}\right)\ldots \left(s_3\ldots s_{l-2}s_{l}s_{l-1}s_{l-2}\ldots s_3\right)\\
   &\qquad\left(s_2\ldots s_{l-2}s_{l}s_{l-1}s_{l-2}\ldots s_2\right)\left(s_1s_2\ldots s_{l-2}s_{l}s_{l-1}s_{l-2}\ldots s_2 s_1\right),
\end{align*}
where the entry $a_{l-1, l}$ in the Cartan matrix equals 0.
When $\g=\nysp$ or $\g=\nyso$,
\begin{align*}
  w_0&=s_{l}\left(s_{l-1}s_ls_{l-1}\right)\left(s_{l-2}s_{l-1}s_{l}s_{l-1}s_{l-2}\right)\ldots \\
   &  \left(s_2\ldots s_{l-1}s_{l}s_{l-1}\ldots s_2\right)\left(s_1s_2\ldots s_{l-1}s_{l}s_{l-1}\ldots s_2 s_1\right),
\end{align*}
where the root $\alpha_l$ is a long root and a short root, respectively.
\end{remark}

When $\g=\nysp$, a path from $v^{+}_1$ to $v^{-}_1$ is given as follows.
\begin{align*}
   v^{-}_1= &\Big(x_{i,0}^{-}\ldots x_{l-1,0}^{-}x_{l,0}^{-}x_{l-1,0}^{-}\ldots x_{i,0}^{-}\Big) \left(\left(x_{i-1,0}^{-}\right)^2x_{i,0}^{-}\ldots x_{l-1,0}^{-}x_{l,0}^{-}x_{l-1,0}^{-}\ldots x_{i,0}^{-}\right)\ldots\\
   &  \left(\left(x_{2,0}^{-}\right)^2\ldots \left(x_{i-1,0}^{-}\right)^2x_{i,0}^{-}\ldots x_{l-1,0}^{-}x_{l,0}^{-}x_{l-1,0}^{-}\ldots x_{i,0}^{-}\right)\\
   &\left(\left(x_{1,0}^{-}\right)^2\ldots \left(x_{i-1,0}^{-}\right)^2x_{i,0}^{-}\ldots x_{l-1,0}^{-}x_{l,0}^{-}x_{l-1,0}^{-}\ldots x_{i,0}^{-}\right)v^{+}_1.
\end{align*}

The degree of the associated polynomial of $Y_{r_j}\left(v_{\sigma_j(\omega_i)}\right)$ is $m_j$. If $m_j=1$, then let $Q_{i,j}(u)=u-a$.
Note that
$\frac{Q_{i,j}\left(u+1\right)}{Q_{i,j}\left(u\right)}= \frac{u-\left(a-1\right)}{u-a}=1+u^{-1}+au^{-2}+a^2u^{-3}+\ldots.$
The eigenvalue of $v_{\sigma_j(\omega_i)}$ under $h_{r_j,1}$ will tell us the value of $a$.
If $m_j=2$, then let $Q_{i,j}(u)=(u-a)(u-b)$.
Note that \begin{align}\label{e5.1}
\frac{Q_{i,j}\left(u+1\right)}{Q_{i,j}\left(u\right)}&=\frac{u-\left(a-1\right)}{u-a}\frac{u-\left(b-1\right)}{u-b} \nonumber\\
&= \left(1+u^{-1}+au^{-2}+a^2u^{-3}+\ldots\right)\left(1+u^{-1}+bu^{-2}+b^2u^{-3}+\ldots\right)\nonumber\\
&= 1+2u^{-1}+\left(a+b+1\right)u^{-2}+\left(a^2+b^2+a+b\right)u^{-3}+\ldots.
\end{align}
Thus the values of both $b$ and $a$ will be derived from the eigenvalues of $v_{\sigma_j(\omega_i)}$ under $h_{r_j,1}$ and $h_{r_j,2}$.


For $i\leq k<l$, denote $x_{k,0}^{-}x_{k+1,0}^{-}\ldots x_{l-1,0}^{-}x_{l,0}^{-}x_{l-1,0}^{-}\ldots x_{i+1,0}^{-}x_{i,0}^{-}$ by $\overline{x_{k,0}^{-}\ldots x_{i,0}^{-}}$ and for $k<i$, denote $\left(x_{k,0}^{-}\right)^2\left(x_{k+1,0}^{-}\right)^2\ldots \left(x_{i-1,0}^{-}\right)^2x_{i,0}^{-}$ $x_{i+1,0}^{-}\ldots x_{l-1,0}^{-}x_{l,0}^{-}x_{l-1,0}^{-}\ldots x_{i+1,0}^{-}x_{i,0}^{-}$ by $\overline{\left(x_{k,0}^{-}\right)^2\ldots x_{i,0}^{-}}$.
In our computations, there is a difference between the cases $1\leq i\leq l-1$ and $i=l$. (Case 2 starts after Proposition \ref{yspc1p111}.)

\noindent\textbf{Case 1:  $1\leq i\leq l-1$.}
\begin{proposition}\label{spc1il-1p1} Let $i\leq k\leq l-1$ and $1\leq m\leq i-2$.\\
\begin{tabular}{|c|c|c|}
  \hline
  Item & $\ysl$-module & Associated polynomial\\\hline
  1 & $Y_{k}\big(x_{k-1,0}^{-}x_{k-2,0}^{-}\ldots x_{i,0}^{-}v^{+}_1\big)$ & $u-\left(a_1+\frac{k-i}{2}\right)$ \\\hline
  2 & $Y_{l}\Big( {x_{l-1,0}^{-}x_{l-2,0}^{-}\ldots x_{i,0}^{-}}v^{+}_1\Big)$ &$u-\frac{1}{2}\big(a_1+\frac{l-i+1}{2}\big)$ \\\hline
  3 & $Y_{k}\left(\overline{x_{k+1,0}^{-}\ldots x_{i,0}^{-}}v^{+}_1\right)$ & $ u-\left(a_1+\frac{2l-i-k+2}{2}\right)$ \\\hline
  4 & $Y_{i-1}\Big(\overline{x_{i,0}^{-}\ldots x_{i,0}^{-}}v^{+}_1\Big)$ & $\Big(u-\left(a_1+l-i+\frac{3}{2}\right)\Big)\Big(u-\left(a_1+\frac{1}{2}\right)\Big)$\\\hline
  5 & $Y_m\Big(\overline{\left(x_{m+1,0}^{-}\right)^2\ldots x_{i,0}^{-}}v^{+}_1\Big)$ & $\Big(u-\left(a_1+\frac{2l-i-m+2}{2}\right)\Big)\Big(u-\left(a_1+\frac{i-m}{2}\right)\Big)$ \\\hline
  6 & $Y_i\Big(\overline{\left(x_{1,0}^{-}\right)^2\ldots x_{i,0}^{-}}v^{+}_1\Big)$ & $u-(a_1+1)$\\
  \hline
\end{tabular}
\end{proposition}
\begin{proof}
The first item is proved in Lemma \ref{ylviw2a-1}. Lemma \ref{ylviw2a0} is devoted to proving the second item, and Lemmas \ref{ylviw2a} and \ref{yspc1il4} handle the third item. The fourth is proved in Lemma \ref{yspc1i5} and the fifth is proved in Lemmas \ref{yspc1i6} and \ref{yspc1i7}. The proof of the last item  is similar to Lemma \ref{yspc1i6}, so we omit the proof.
\end{proof}

\begin{lemma}\label{ylviw2a-1}
The associated polynomial of $Y_{k}\Big(\big(x_{k-1,0}^{-}x_{k-2,0}^{-}\ldots x_{i,0}^{-}\big)v^{+}_1\Big)$ is given by $u-\left(a_1+\frac{k-i}{2}\right)$ for $i\leq k\leq l-1$.
\end{lemma}
\begin{proof}
The associated polynomial of $Y_{k}\Big(\big(x_{k-1,0}^{-}x_{k-2,0}^{-}\ldots x_{i,0}^{-}\big)v^{+}_1\Big)$  is of degree 1, say $u-a_k$. The value $a_k$ is the eigenvalue of $x_{k-1,0}^{-}\ldots x_{i,0}^{-}v^{+}_1$ under $h_{k,1}$. We claim $a_k=a_1+\frac{k-i}{2}$ and prove this by using induction on $k$. If $k=i$, then $h_{i,1}v^{+}_1=a_1v^{+}_1$ and hence the claim is true. Suppose this is true for $k-1$.
By the induction hypothesis, we have $$h_{k-1,1}x_{k-2,0}^{-}\ldots x_{i,0}^{-}v^{+}_1=\left(a_1+\frac{k-1-i}{2}\right)x_{k-2,0}^{-}\ldots x_{1,0}^{-}v^{+}_1,$$ then by $(\ref{Cor3.3})$ $$x_{k-1,1}^{-}x_{k-2,0}^{-}\ldots x_{i,0}^{-}v^{+}_1=\left(a_1+\frac{k-1-i}{2}\right)x_{k-1,0}^{-}x_{k-2,0}^{-}\ldots x_{i,0}^{-}v^{+}_1.$$
To show that the claim is true for $k$, we note
\begin{align*}
   h_{k,1}x_{k-1,0}^{-}&x_{k-2,0}^{-}\ldots x_{1,0}^{-}v^{+}_1 \\
   &= [h_{k,1},x_{k-1,0}^{-}]x_{k-2,0}^{-}\ldots x_{1,0}^{-}v^{+}_1\\
 &= \left(x_{k-1,1}^{-}+\frac{1}{2}x_{k-1,0}^{-}+x_{k-1,0}^{-}h_{k,0}\right)x_{k-2,0}^{-}\ldots x_{1,0}^{-}v^{+}_1\\
 &= \left(a_1+\frac{k-1-i}{2}+\frac{1}{2}\right)x_{k-1,0}^{-}x_{k-2,0}^{-}\ldots x_{1,0}^{-}v^{+}_1\\
 &= \left(a_1+\frac{k-i}{2}\right)x_{k-1,0}^{-}x_{k-2,0}^{-}\ldots x_{1,0}^{-}v^{+}_1.
\end{align*}
Therefore the claim is proved by induction.
\end{proof}
Recall that $\{x_{l,r}^{\pm}, h_{l,r}|r\in\Z_{\geq 0}\}$ generates a subalgebra $Y_l$ of $\yg$ and $Y_l\cong \ysl$. However, $x_{l,r}^{\pm}$ and $h_{l,r}$ do not satisfy the defining relations of $\ysl$. Therefore, in order to apply the results in the case of $\ysl$,  we need to rescale these generators.
Let $\tilde{x}_{l,r}^{\pm}=\frac{\sqrt{2}}{2^{r+1}}x_{l,r}^{\pm}$ and $\tilde{h}_{l,r}=\frac{1}{2^{r+1}}h_{l,r}.$ The new generators satisfy the defining relations of $\ysl$.

\begin{lemma}\label{ylviw2a0}
The associated polynomial of $Y_{l}\big( {x_{l-1,0}^{-}\ldots x_{i+1,0}^{-}x_{i,0}^{-}}v^{+}_1\big)$ is given by $u-\frac{1}{2}\left(a_1+\frac{l-i+1}{2}\right)$.
\end{lemma}
\begin{proof}The associated polynomial  is of degree 1, say $u-a$. The constant $a$ is the eigenvalue of $x_{l-1,0}^{-}x_{l-2,0}^{-}\ldots x_{i,0}^{-}v^{+}_1$ under $\widetilde{h}_{l,1}$.
\begin{align*}
  \tilde{h}_{l,1}\big(x_{l-1,0}^{-}&x_{l-2,0}^{-}\ldots x_{i,0}^{-}\big)v^{+}_1\\
 &= \frac{1}{4}h_{l,1}\left(x_{l-1,0}^{-}x_{l-2,0}^{-}\ldots x_{i,0}^{-}\right)v^{+}_1\\
 &= \frac{1}{4}[h_{l,1},x_{l-1,0}^{-}]\left(x_{l-2,0}^{-}\ldots x_{i,0}^{-}\right)v^{+}_1\\
 &=  \frac{1}{4}\left(2x_{l-1,1}^{-}+h_{l,0}x_{l-1,0}^{-}+x_{l-1,0}^{-}h_{l,0}\right)x_{l-2,0}^{-}\ldots x_{i,0}^{-}v^{+}_1 \\
 &=  \frac{1}{4}\left(2x_{l-1,1}^{-}+h_{l,0}x_{l-1,0}^{-}\right)x_{l-2,0}^{-}\ldots x_{i,0}^{-}v^{+}_1 \\
 &=  \frac{1}{2}\left(a_1+\frac{l-i-1}{2}+1\right)x_{l-1,0}^{-}x_{l-2,0}^{-}\ldots x_{i,0}^{-}v^{+}_1\quad \big(\text{by (4.1)}\big)\\
  &=  \frac{1}{2}\left(a_1+\frac{l-i+1}{2}\right)x_{l-1,0}^{-}x_{l-2,0}^{-}\ldots x_{i,0}^{-}v^{+}_1.
\end{align*}
\end{proof}

\begin{lemma}\label{ylviw2a} The associated polynomial of
$Y_{l-1}\left(x_{l,0}^{-}x_{l-1,0}^{-}\ldots x_{i+1,0}^{-}x_{i,0}^{-}v^{+}_1\right)$ is given by $ u-\left(a_1+\frac{l-i+3}{2}\right)$.
\end{lemma}

\begin{proof}The associated polynomial is of degree 1, say $u-a$. The value $a$ equals the eigenvalue of $x_{l,0}^{-}x_{l-1,0}^{-}\ldots x_{i,0}^{-}v^{+}_1$ under $h_{l-1,1}$. By Lemma \ref{ylviw2a0} and (4.1), we have that $x_{l,1}^{-}x_{l-1,0}^{-}\ldots x_{i,0}^{-}v^{+}_1=2\sqrt{2}\tilde{x}_{l,1}^{-}\ldots x_{i,0}^{-}v^{+}_1=\sqrt{2}\left(a_1+\frac{l-i+1}{2}\right)\widetilde{x}_{l,0}^{-}\ldots x_{i,0}^{-}v^{+}_1=\left(a_1+\frac{l-i+1}{2}\right)x_{l,0}^{-}\ldots x_{i,0}^{-}v^{+}_1.$
\begin{align*}
   h&_{l-1,1}\big(x_{l,0}^{-}x_{l-1,0}^{-}\ldots x_{i,0}^{-}\big)v^{+}_1 \\
   &= [h_{l-1,1}x_{l,0}^{-}]x_{l-1,0}^{-}\ldots x_{i,0}^{-}v^{+}_1+x_{l,0}^{-}h_{l-1,1}x_{l-1,0}^{-}\ldots x_{i,0}^{-}v^{+}_1 \\
   &= \left(2x_{l,1}^{-}+2x_{l,0}^{-}+2x_{l,0}^{-}h_{l-1,0}-\left(a_1+\frac{l-i-1}{2}\right)x_{l,0}^{-}\right)x_{l-1,0}^{-}\ldots x_{i,0}^{-}v^{+}_1 \quad \big(\text{by (4.1)}\big)\\
   &= 2x_{l,1}^{-}x_{l-1,0}^{-}\ldots x_{i,0}^{-}v^{+}_1-\left(a_1+\frac{l-i-1}{2}\right)x_{l,0}^{-}x_{l-1,0}^{-}\ldots x_{i,0}^{-}v^{+}_1\\
   &= 2\left(a_1+\frac{l-i+1}{2}\right)x_{l,0}^{-}x_{l-1,0}^{-}\ldots x_{i,0}^{-}v^{+}_1 -\left(a_1+\frac{l-i-1}{2}\right)x_{l,0}^{-}x_{l-1,0}^{-}\ldots x_{i,0}^{-}v^{+}_1\\
   &=\left(a_1+\frac{l-i+3}{2}\right)x_{l,0}^{-}x_{l-1,0}^{-}\ldots x_{i,0}^{-}v^{+}_1.
\end{align*}
\end{proof}

\begin{lemma}\label{yspc1il4}
For $i\leq k\leq l-2$, the associated polynomial of $Y_{k}\left(\overline{x_{k+1,0}^{-}\ldots x_{i,0}^{-}}v^{+}_1\right)$ is $ u-\left(a_1+\frac{2l-i-k+2}{2}\right)$.
\end{lemma}
\begin{proof} The proof is similar to the proof of Lemma \ref{ylviw2a}, so we omit the proof.
\end{proof}
If $i=1$, the proof of Proposition \ref{spc1il-1p1} is complete and we move on to Case 2. Now suppose that $i\geq 2$.

\begin{lemma}\label{yspc1i5}
The associated polynomial of the representation $Y_{i-1}\left(\overline{x_{i,0}^{-}\ldots x_{i,0}^{-}}v^{+}_1\right)$ is given by $ \Big(u-\left(a_1+l-i+\frac{3}{2}\right)\Big)\Big(u-\left(a_1+\frac{1}{2}\right)\Big)$.
\end{lemma}
\begin{proof}
The associated polynomial of the representation $Y_{i-1}\left(\overline{x_{i,0}^{-}\ldots x_{i,0}^{-}}v^{+}_1\right)$ has degree 2, say $\left(u-a\right)\left(u-b\right)$.
The eigenvalues of $\overline{x_{i,0}^{-}\ldots x_{i,0}^{-}}v^{+}_1$ under $h_{i-1,1}$ and $h_{i-1,2}$ will tell us the values of $a$ and $b$ (see (\ref{e5.1})).
Let $A=a_1+l-i+1$ and $B=a_1$.
\begin{align*}
   h_{i-1,1}&\overline{x_{i,0}^{-}\ldots x_{i,0}^{-}}v^{+}_1\\
   &= [h_{i-1,1},x_{i,0}^{-}]\overline{x_{i+1,0}^{-}\ldots x_{i,0}^{-}}v^{+}_1 +\overline{x_{i,0}^{-}\ldots x_{i+1,0}^{-}}[h_{i-1,1},x_{i,0}^{-}]v^{+}_1\\
    &= \left(x_{i,1}^{-}+\frac{1}{2}x_{i,0}^{-}+x_{i,0}^{-}h_{i-1,0}\right)\overline{x_{i+1,0}^{-}\ldots x_{i,0}^{-}}v^{+}_1\\
     &+\overline{x_{i,0}^{-}\ldots x_{i+1,0}^{-}}\left(x_{i,1}^{-}+\frac{1}{2}x_{i,0}^{-}+x_{i,0}^{-}h_{i-1,0}\right)v^{+}_1 \\
     &= \Bigg(\left(a_1+\frac{2l-i-i+2}{2}+\frac{1}{2}+1\right)+\left(a_1+\frac{1}{2}\right)\Bigg)\overline{x_{i,0}^{-}\ldots x_{i,0}^{-}}v^{+}_1\\
      &= \left(\left(A+\frac{1}{2}\right)+\left(B+\frac{1}{2}\right)+1\right)\overline{x_{i,0}^{-}\ldots x_{i,0}^{-}}v^{+}_1.
\end{align*}
To obtain the third equality, we used Lemma \ref{yspc1il4} when $k=i$ and (4.1).
\begin{align*}
h&_{i-1,2}\overline{x_{i,0}^{-}\ldots x_{i,0}^{-}}v^{+}_1 \\
&= [h_{i-1,2},x_{i,0}^{-}]\overline{x_{i+1,0}^{-}\ldots x_{i,0}^{-}}v^{+}_1+\overline{x_{i,0}^{-}\ldots x_{i+1,0}^{-}}[h_{i-1,2},x_{i,0}^{-}]v^{+}_1\\
&= \left(x_{i,2}^{-}+x_{i,1}^{-}+\frac{1}{4}x_{i,0}^{-}+x_{i,1}^{-}h_{i-1,0}+
\frac{1}{2}x_{i,0}^{-}h_{i-1,0}+x_{i,0}^{-}h_{i-1,1}\right)\overline{x_{i+1,0}^{-}\ldots x_{i,0}^{-}}v^{+}_1 \\
&+\overline{x_{i,0}^{-}\ldots x_{i+1,0}^{-}}\left(x_{i,2}^{-}+x_{i,1}^{-}+\frac{1}{4}x_{i,0}^{-}+x_{i,1}^{-}h_{i-1,0}+
\frac{1}{2}x_{i,0}^{-}h_{i-1,0}+x_{i,0}^{-}h_{i-1,1}\right)v^{+}_1 \\
&= \left(x_{i,2}^{-}+2x_{i,1}^{-}+\frac{3}{4}x_{i,0}^{-}\right)\overline{x_{i+1,0}^{-}\ldots x_{i,0}^{-}}v^{+}_1+\overline{x_{i,0}^{-}\ldots x_{i+1,0}^{-}}h_{i-1,1}x_{i,0}^{-}v^{+}_1\\
&+\overline{x_{i,0}^{-}\ldots x_{i+1,0}^{-}}\left(x_{i,2}^{-}+x_{i,1}^{-}+\frac{1}{4}x_{i,0}^{-}\right)v^{+}_1 \\
&= \left(x_{i,2}^{-}+2x_{i,1}^{-}+\frac{3}{4}x_{i,0}^{-}\right)\overline{x_{i+1,0}^{-}\ldots x_{i,0}^{-}}v^{+}_1+\overline{x_{i,0}^{-}\ldots x_{i+1,0}^{-}}\left(x_{i,2}^{-}+2x_{i,1}^{-}+\frac{3}{4}x_{i,0}^{-}\right)v^{+}_1 \\
&=\Big(A^2+2A+B^2+2B+\frac{3}{2}\Big)\overline{x_{i,0}^{-}\ldots x_{i,0}^{-}}v^{+}_1\\
&=\left(\left(A+\frac{1}{2}\right)^2+\left(B+\frac{1}{2}\right)^2+\left(A
+\frac{1}{2}\right)+\left(B+\frac{1}{2}\right)\right)\overline{x_{i,0}^{-}\ldots x_{i,0}^{-}}v^{+}_1.
\end{align*}
Therefore it follows from (5.1) that $a=a_1+\frac{1}{2}$ and $b=a_1+l-i+\frac{3}{2}$, or vice-versa with $a$ and $b$ switched.
\end{proof}


For $i=2$, see Lemma \ref{sps1v2p}. Now we suppose $i\geq 3$.
\begin{lemma}\label{yspc1i6}
The associated polynomial $Q(u)$ of $Y_{i-2}\Big(\overline{\left(x_{i-1,0}^{-}\right)^2\ldots x_{i,0}^{-}}v^{+}_1\Big)$ is given by $\Big(u-\left(a_1+l-i+2\right)\Big)\Big(u-\left(a_1+1\right)\Big)$.
\end{lemma}
\begin{proof} 
The associated polynomial $Q\left(u\right)$ has degree 2. Let $Q\left(u\right)=\left(u-a\right)\left(u-b\right)$.
The eigenvalues of $\overline{\left(x_{i-1,0}^{-}\right)^2\ldots x_{i,0}^{-}}v^{+}_1$ under $h_{i-2,1}$ and $h_{i-2,2}$ will tell us the values of $a$ and $b$. Denote $a_1+l-i+\frac{3}{2}$ by $A$ and $a_1+\frac{1}{2}$ by $B$.  By the defining relations, $h_{i-2,0}\overline{x_{i-1,0}^{-}\ldots x_{i,0}^{-}}v^{+}_1=\overline{x_{i-1,0}^{-}\ldots x_{i,0}^{-}}v^{+}_1$.
\begin{align*}
h_{i-2,1}&\overline{\left(x_{i-1,0}^{-}\right)^2\ldots x_{i,0}^{-}}v^{+}_1 \\
   &= \left[h_{i-2,1}, \left(x_{i-1,0}^{-}\right)^2\right]\overline{x_{i,0}^{-}\ldots x_{i,0}^{-}}v^{+}_1\\
   &= [h_{i-2,1}, x_{i-1,0}^{-}]\overline{x_{i-1,0}^{-}\ldots x_{i,0}^{-}}v^{+}_1+ x_{i-1,0}^{-}[h_{i-2,1}, x_{i-1,0}^{-}]\overline{x_{i,0}^{-}\ldots x_{i,0}^{-}}v^{+}_1\\
   &=\left(x_{i-1,1}^{-}+\frac{1}{2}x_{i-1,0}^{-}+x_{i-1,0}^{-}h_{i-2,0}\right)\overline{x_{i-1,0}^{-}\ldots x_{i,0}^{-}}v^{+}_1\\
   &+x_{i-1,0}^{-}\left(x_{i-1,1}^{-}+\frac{1}{2}x_{i-1,0}^{-}+x_{i-1,0}^{-}h_{i-2,0}\right)\overline{x_{i,0}^{-}\ldots x_{i,0}^{-}}v^{+}_1\\
   &= \left(x_{i-1,1}^{-}x_{i-1,0}^{-}+x_{i-1,0}^{-}x_{i-1,1}^{-}\right)\overline{x_{i,0}^{-}\ldots x_{i,0}^{-}}v^{+}_1+2\overline{\left(x_{i-1,0}^{-}\right)^2\ldots x_{i,0}^{-}}v^{+}_1\\
   &= \left(A+B+2\right)\overline{\left(x_{i-1,0}^{-}\right)^2\ldots x_{i,0}^{-}}v^{+}_1\\
    &= \left(\left(A+\frac{1}{2}\right)+\left(B+\frac{1}{2}\right)+1\right)\overline{\left(x_{i-1,0}^{-}\right)^2\ldots x_{i,0}^{-}}v^{+}_1,
\end{align*}
where the second equality from the end follows from Corollary \ref{w1bw1a} and Lemma \ref{yspc1i5}.

\begin{align*}
   h_{i-2,2}&\overline{\left(x_{i-1,0}^{-}\right)^2\ldots x_{i,0}^{-}}v^{+}_1 \\
   &= \left[h_{i-2,2}, \left(x_{i-1,0}^{-}\right)^2\right]\overline{x_{i,0}^{-}\ldots x_{i,0}^{-}}v^{+}_1 \\
   &= [h_{i-2,2}, x_{i-1,0}^{-}]\overline{x_{i-1,0}^{-}\ldots x_{i,0}^{-}}v^{+}_1+ x_{i-1,0}^{-}[h_{i-2,2}, x_{i-1,0}^{-}]\overline{x_{i,0}^{-}\ldots x_{i,0}^{-}}v^{+}_1\\
   &=\left(x_{i-1,2}^{-}x_{i-1,0}^{-}+x_{i-1,0}^{-}x_{i-1,2}^{-}+2x_{i-1,1}^{-}x_{i-1,0}^{-}+2x_{i-1,0}^{-}x_{i-1,1}^{-}
   +\frac{3}{2}\left(x_{i-1,0}^{-}\right)^2\right)\\
   &\qquad\qquad\overline{x_{i,0}^{-}\ldots x_{i,0}^{-}}v^{+}_1\\
   &=\Big(A^2+B^2+2\left(A+B\right)+\frac{3}{2}\Big)\overline{\left(x_{i-1,0}^{-}\right)^2\ldots x_{i,0}^{-}}v^{+}_1\\
   &=\Bigg(\left(A+\frac{1}{2}\right)^2+\left(B+\frac{1}{2}\right)^2+\left(A
+\frac{1}{2}\right)+\left(B+\frac{1}{2}\right)\Bigg)\overline{\left(x_{i-1,0}^{-}\right)^2\ldots x_{i,0}^{-}}v^{+}_1,
\end{align*}
where the second equality from the end follows from Corollary \ref{w1bw1a} and Lemma \ref{yspc1i5}. Therefore it follows from (5.1) that  $a=a_1+1$ and $b=a_1+l-i+2$, or vice-versa with $a$ and $b$ switched.
\end{proof}

Similarly to Lemma \ref{yspc1i6}, using induction on $m$ downward, we have

\begin{lemma}\label{yspc1i7}
If $1\leq m\leq i-3$, the associated polynomial of
$Y_m\big(\overline{(x_{m+1,0}^{-})^2\ldots x_{i,0}^{-}}v^{+}_1\big)$ is given by $\big(u-(a_1+\frac{2l-i-m+2}{2})\big)\big(u-(a_1+\frac{i-m}{2})\big)$.
\end{lemma}
The following lemma follows immediately from explicit computations.
\begin{lemma}\label{sps1v2p}
Let $1\leq n\leq i-1$. Denote $s_ns_{n+1}\ldots s_{l-1}s_ls_{l-1}\ldots s_i$ by $\mathbf{s}_n$. Then $\mathbf{s}_n^{-1}\left(\alpha_j\right)$ is a positive root, where $\alpha_j$ is a simple root of $\nysp$ and $j=n+1,n+2,\ldots, l$.
\end{lemma}

Using the idea of the next proposition, we will be able to obtain also Proposition \ref{yspc1p111}, and in this way we will not have to repeat the previous long computations needed to prove Proposition 5.5.
Let $v_2=\overline{(x_{1,0}^{-})^2\ldots x_{i,0}^{-}}v^{+}_1$.

\begin{proposition}\label{spc1il-1p2} Let $i\leq k\leq l-1$ and $2\leq m\leq i-2$.\\
\begin{tabular}{|c|c|c|}
  \hline
  Item & $\ysl$-module & Associated polynomial\\\hline
  1 & $Y_{k}\big(x_{k-1,0}^{-}x_{k-2,0}^{-}\ldots x_{i,0}^{-}v_2\big)$ & $u-\left(a_1+1+\frac{k-i}{2}\right)$ \\\hline
  2 & $Y_{l}\Big( {x_{l-1,0}^{-}x_{l-2,0}^{-}\ldots x_{i,0}^{-}}v_2\Big)$ &$u-\frac{1}{2}\big(a_1+1+\frac{l-i+1}{2}\big)$ \\\hline
  3 & $Y_{k}\left(\overline{x_{k+1,0}^{-}\ldots x_{i,0}^{-}}v_2\right)$ & $ u-\left(a_1+1+\frac{2l-i-k+2}{2}\right)$ \\\hline
  4 & $Y_{i-1}\Big(\overline{x_{i,0}^{-}\ldots x_{i,0}^{-}}v_2\Big)$ & $\Big(u-\left(a_1+l-i+\frac{5}{2}\right)\Big)\Big(u-\left(a_1+1+\frac{1}{2}\right)\Big)$\\\hline
  5 & $Y_m\Big(\overline{\left(x_{m+1,0}^{-}\right)^2\ldots x_{i,0}^{-}}v_2\Big)$ & $\Big(u-\left(a_1+\frac{2l-i-m+4}{2}\right)\Big)\Big(u-\left(a_1+1+\frac{i-m}{2}\right)\Big)$ \\\hline
  6 & $Y_i\Big(\overline{\left(x_{1,0}^{-}\right)^2\ldots x_{i,0}^{-}}v_2\Big)$ & $u-(a_1+2)$\\
  \hline
\end{tabular}
\end{proposition}
\begin{proof}
By removing the first node in the Dynkin diagram of the Lie algebra of type $C_l$, we obtain a simple Lie algebra of type $C_{l-1}$. Let $I'=\{2,3,\ldots, l\}$.  Let $Y^{\left(1\right)}$ be the Yangian generated by all $x_{j,k}^{\pm}$ and $h_{j,k}$, for $j\in I'$ and $k\in \mathbb{Z}_{\geq 0}$, so that $Y^{(1)}\cong Y\big(\mathfrak{sp}(2(l-1),\C)\big)$.

As in Lemma \ref{S5L1yiahwr}, we can show, by Lemma \ref{sps1v2p}, that $Y^{\left(1\right)}\left(v_2\right)$ is a highest weight representation of $Y\big(\mathfrak{sp}(2(l-1),\C)\big)$. Suppose that the associated $(l-1)$-tuple of polynomials is $\pi'=\big(\pi'_2(u),\ldots, \pi'_l(u)\big)$. Since the weight of $v_2$ is $-2\omega_1+\omega_i$, $\pi'_j(u)=1$ if $j\neq i$ and $\pi'_i(u)=u-a$. It follows from the sixth item of Proposition \ref{spc1il-1p1} that $a=a_1+1$. The remainder of the proof of this proposition can be obtained by replacing $a_1$ with $a_1+1$ in the proof of Proposition \ref{spc1il-1p1}. This explains why this proposition is the same as Proposition \ref{spc1il-1p1}, but with $a_1$ replaced by $a_1+1$.
\end{proof}

Note that if $i=2$, only the first three items in the above proposition are necessary.

We now assume that $i\geq 3$.
Inductively define $v_{n+1}=\overline{\left(x_{n,0}^{-}\right)^2\ldots x_{i,0}^{-}}v_{n}$ for $2\leq n\leq i-1$, and let $Y^{\left(n\right)}$ be the Yangian generated by all $x_{r,k}^{\pm}$ and $h_{r,k}$ for $r>n$ and $k\in \mathbb{Z}_{\geq 0}$. Since $n\leq i-1\leq l-2$, $Y^{(n)}$ is a Yangian of type $C_{l-n}$. We can show
\begin{proposition}\label{yspc1p111} Let $i\leq k\leq l-1$ and $n+1\leq m\leq i-2$.\\
{\small \begin{tabular}{|c|c|c|}
  \hline
  item & $\ysl$-module & Associated polynomial\\\hline
  1 & $Y_{k}\big(x_{k-1,0}^{-}x_{k-2,0}^{-}\ldots x_{i,0}^{-}v_n\big)$ & $u-\left(a_1+n+\frac{k-i}{2}\right)$ \\\hline
  2 & $Y_{l}\Big( {x_{l-1,0}^{-}x_{l-2,0}^{-}\ldots x_{i,0}^{-}}v_n\Big)$ &$u-\frac{1}{2}\big(a_1+n+\frac{l-i+1}{2}\big)$ \\\hline
  3 & $Y_{k}\left(\overline{x_{k+1,0}^{-}\ldots x_{i,0}^{-}}v_n\right)$ & $ u-\left(a_1+n+\frac{2l-i-k+2}{2}\right)$ \\\hline
  4 & $Y_{i-1}\Big(\overline{x_{i,0}^{-}\ldots x_{i,0}^{-}}v_n\Big)$ & $\Big(u-\left(a_1+n+l-i+\frac{3}{2}\right)\Big)\Big(u-\left(a_1+n+\frac{1}{2}\right)\Big)$\\\hline
  5 & $Y_m\Big(\overline{\left(x_{m+1,0}^{-}\right)^2\ldots x_{i,0}^{-}}v_n\Big)$ & $\Big(u-(a_1+n+\frac{2l-i-m+2}{2})\Big)\Big(u-(a_1+n+\frac{i-m}{2})\Big)$ \\\hline
  6 & $Y_i\Big(\overline{\left(x_{1,0}^{-}\right)^2\ldots x_{i,0}^{-}}v_n\Big)$ & $u-(a_1+n+1)$\\
  \hline
\end{tabular}}
\end{proposition}


\noindent\textbf{ Case 2: $i=l$.}\\
Denote $v_{1}^{+}$ by $v_{1}$ and inductively define $v_{n+1}=\left(x_{n,0}^{-}\right)^2\left(x_{n+1,0}^{-}\right)^2\ldots \left(x_{l-1,0}^{-}\right)^2x_{l,0}^{-}v_{n}$ for $1\leq n\leq l-2$. We omit the proof of the following proposition because it is very similar to the proof of Lemmas \ref{yspc1i5} and \ref{yspc1i6} and Proposition \ref{spc1il-1p2}.
\begin{proposition}\label{yspc3p1} Let $m=l-1$ or $1\leq m\leq l-2$ and $m\leq k\leq l-2$.\\
\begin{tabular}{|c|c|c|}
  \hline
  Item & $\ysl$-module & Associated polynomial\\\hline
  1 & $Y_l\left(v_{m}\right)$ & $u-\frac{a_1+m-1}{2}$ \\\hline
  2 & $Y_{l-1}\left(x_{l,0}^{-}v_m\right)$ &$\Big(u-(a_1+m)\Big)\Big(u-(a_1+m-1)\Big)$ \\\hline
  3 & $Y_{k}\left(\overline{(x_{k+1,0}^{-})^2\ldots x_{l,0}^{-}}v_m\right)$ & $\Big(u-(a_1+m+\frac{l-k-1}{2})\Big)\Big(u-(a+m+\frac{l-k-3}{2})\Big)$\\\hline
\end{tabular}
\end{proposition}
We summarize all results in Cases 1 and 2 into the following corollary.
\begin{corollary}\label{aprcsp12} The set $T(i, r_j)$ of all possible roots of the associated polynomial of $Y_{r_j}\left(v_{\sigma_j(\omega_i)}\right)$ is listed below. For $1\leq i, r_j\leq l-1$,
\begin{enumerate}
  \item $T\left(i,r_j\right)=\left\{a_1+\frac{|i-r_j|}{2}+r, a_1+l+1+r-\frac{i+r_j}{2}|0\leq r<\text{min}\{i, r_j\}\right\}$;

  \item $T\left(i,l\right)=\left\{\frac{1}{2}\left(a_1+\frac{l-i+1}{2}+r\right)|0\leq r< i\right\};$

  \item  $T\left(l,r_j\right)=\left\{a_1+\frac{l-r_j+1}{2}+r, a_1+\frac{l-r_j-1}{2}+r|0\leq r<r_j\right\}$;

   \item $T\left(l,l\right)=\left\{\frac{a_1}{2}, \frac{a_1+1}{2},\ldots, \frac{a_1+l-1}{2} \right\}$.
\end{enumerate}
\end{corollary}
Similar computations can be carried out  when $\g$ is a classical simple Lie algebra of any other type to obtain a concrete set $T(i, r_j)$. We omit the details here, see \cite{Ta}.
By Theorem \ref{S5T1Liahw}, using the sets $T(i, r_j)$, we obtain a concrete cyclicity condition for certain tensor products.
\begin{theorem}\label{mtoyspl}
Let $L=V_{a_1}(\omega_{b_1})\otimes V_{a_2}(\omega_{b_2})\otimes\ldots\otimes V_{a_k}(\omega_{b_k})$ be an ordered tensor product of fundamental representations of $\yg$. If $a_n-a_m\notin S\left(b_m,b_n\right)$ for $n>m$, then $L$ is a highest weight representation of $\yg$, where the set $S\left(b_m,b_n\right)$ is defined as follows:
\begin{enumerate}

  \item When $\g=\nyn$,
$S\left(b_m,b_n\right)=\left\{\frac{|b_n-b_m|}{2}+k|1\leq k\leq \text{min}\left\{b_m, l-b_n+1\right\}\right\}$.

\item When $\g=\nyso$,
  \begin{enumerate}
  \item $S\left(b_m,b_n\right)=\left\{|b_m-b_n|+2+2r, 2l-(b_m+b_n)+1+2r\right\}$, \\
   where $1\leq b_m, b_n\leq l-1$ and $0\leq r<\text{min}\{b_m, b_n\}$;
  \item  $S\left(l,b_n\right)=\left\{l-b_n+2+2r|0\leq r<b_n, 1\leq b_n\leq l-1\right\}$;
  \item $S\left(b_m,l\right)=\left\{l-b_m+1+r, l-b_m+r|0\leq r< b_m, 1\leq b_m\leq l-1\right\}$;
   \item $S\left(l,l\right)=\left\{1,3,\ldots,2l-1\right\}$.
\end{enumerate}

 \item When $\g=\nysp$, \begin{enumerate}
  \item $S\left(b_m,b_n\right)=\left\{\frac{|b_m-b_n|}{2}+1+r, l+2+r-\frac{b_m+b_n}{2}|0\leq r<\text{min}\{b_m, b_n\}\right\}$, \\
   where $1\leq b_m, b_n\leq l-1$;
  \item  $S\left(l,b_n\right)=\left\{\frac{l-b_n+1}{2}+1+r, \frac{l-b_n-1}{2}+1+r|0\leq r<b_n, 1\leq b_n\leq l-1\right\}$;
  \item $S\left(b_m,l\right)=\left\{\frac{l-b_m+1}{2}+2+r|0\leq r< b_m, 1\leq b_m\leq l-1\right\};$
   \item $S\left(l,l\right)=\left\{2,3,\ldots,l+1\right\}$.
\end{enumerate}

  \item When $\g=\nyo$,
  \begin{enumerate}
  \item $S\left(b_m,b_n\right)=\{\frac{|b_m-b_n|}{2}+1+r, l+r-\frac{b_m+b_n}{2}\}$, where $1\leq b_m, b_n\leq l-2$ and $0\leq r<\text{min}\{b_m, b_n\}$;
  \item  $S\left(l-1,b_n\right)=S\left(l,b_n\right)=S\left(b_n,l-1\right)=S\left(b_n,l\right)=\\ \left\{\frac{l-1-b_n}{2}+1+r|0\leq r<b_n, 1\leq b_n\leq l-2\right\}$;
  \item $S\left(l-1,l\right)=S\left(l,l-1\right)=\left\{2,4,\ldots,l-2+\overline{l}\right\}$,  where $\overline{l}=0$ if $l$ is even and let $\overline{l}=1$ if $l$ is odd;
   \item $S\left(l-1,l-1\right)=S\left(l,l\right)=\left\{1,3,\ldots,l-1-\overline{l}\right\}$.
\end{enumerate}
\end{enumerate}
\end{theorem}
\begin{remark}\
\begin{enumerate}
  \item Our result on the cyclicity condition for $L$ is an analogue of the case of Corollary 6.2 in \cite{Ch3} in which $m_1=m_2=1$. The numbers in the sets $S(i_1, i_2)$ as defined above are the exponents of $q$ in the set $\mathcal{S}(i_1, i_2)$ as defined in \cite{Ch3} up to a factor of 2 and a constant.
  \item When $\g=\nyn$ and $k=2$, the previous theorem gives the same condition for cyclicity as Theorem 6.2 in \cite{ChPr8}. when $\g=\nyo$ and $k=2$, our results in cases (b), (c) and (d) above agree with those in Theorem 7.2 in \cite{ChPr8}. In case (a) with $k=2$, if $b_m+b_n\leq l$, our set $S(b_m,b_n)$ is the same as the corresponding set of values in Theorem 7.2 in \cite{ChPr8}; if $b_m+b_n>l$, our set $S(b_m,b_n)$ contains it strictly.
\end{enumerate}
\end{remark}
By Proposition \ref{VoWWoVhi}, $L$ is irreducible if and only if both $L$ and the left dual$\ ^t L$ are highest weight representations. Denote the half of dual Coxeter number of $\g$ by $\kappa$. By Proposition 3.4 and Corollary 3.6 in \cite{ChPr8} that $\ ^t L\cong V_{a_{k}-\kappa}\left(\omega_{-w_0(b_{k})}\right)\otimes V_{a_{k-1}-\kappa}\left(\omega_{-w_0(b_{k-1})}\right)\otimes \ldots \otimes V_{a_{1}-\kappa}\left(\omega_{-w_0(b_{1})}\right)$. The cyclicity condition for $L$ leads to an irreducibility criterion for $L$. We note that when $\g=\nyn$, $S(b_m,b_n)=S(l-b_n+1, l-b_m+1)$ for $1\leq b_m, b_n\leq l$.
\begin{theorem}\label{icotpofr}
Let $L$ and $S(b_i, b_j)$ be defined as in Theorem \ref{mtoyspl}.
If $a_j-a_i\notin S(b_i, b_j)$, for $1\leq i\neq j\leq k$, then $L$ is an irreducible representation of $\yg$.
\end{theorem}
When $\g$ is of type $A$, the implication in Theorem \ref{icotpofr} becomes an equivalence. This result is given in Theorem \ref{taiaeaoic} , but to prove this, we use the following lemma.
\begin{lemma}[Theorem 6.2, \cite{ChPr8}]\label{cp8c3t}
Let $1\leq b_i\leq b_j\leq l$. $V_{a_i}(\omega_{b_i})\otimes V_{a_j}(\omega_{b_j})$ is reducible as a $\yn$-module if and only if
\begin{center}
 $a_j-a_i=\pm\big(\frac{b_j-b_i}{2}+r\big)$, where $0<r\leq\text{min}(b_i,l+1-b_j).$
\end{center}
\end{lemma}
\begin{theorem}\label{taiaeaoic}
Let $L=V_{a_1}(\omega_{b_1})\otimes V_{a_2}(\omega_{b_2})\otimes\ldots\otimes V_{a_k}(\omega_{b_k})$ be an ordered tensor product of fundamental representations of $\yn$.  $L$ is irreducible if and only if $a_j-a_i\notin S\left(b_i,b_j\right)$ for any $a_i$ and $a_j$ with $1\leq i\neq j\leq k$.
\end{theorem}
\begin{proof}
The sufficiency of this condition is proved in Theorem \ref{icotpofr}. We now show the necessity of this condition. Suppose that $L$ is irreducible. Suppose on the contrary that there exist $a_i$ and $a_j$ for $i\neq j$ such that $a_j-a_i\in S(b_i,b_j)$. Since $L$ is irreducible, any permutation of tensor factors $V_{a_s}(\omega_{b_s})$ gives an isomorphic representation of $\yn$. Arranging the order if necessary, we may assume that $a_2-a_1\in S(b_1, b_2)$. Suppose $b_1\leq b_2$. By Lemma \ref{cp8c3t},  $V_{a_1}(\omega_{b_1})\otimes V_{a_2}(\omega_{b_2})$ is reducible, hence so is $L$, which is a contradiction. So $b_1>b_2$.  Note that $S(b_1,b_2)=S(b_2,b_1)$. Thus $a_2-a_1\in S(b_2, b_1)$, i.e., $a_2-a_1=\frac{b_1-b_2}{2}+r$ for some $0<r\leq\text{min}(b_2,l+1-b_1).$ Hence we have $a_1-a_2=-\big(\frac{b_1-b_2}{2}+r\big)$, and then $V_{a_2}(\omega_{b_2})\otimes V_{a_1}(\omega_{b_1})$ is reducible by Lemma \ref{cp8c3t}, which implies that $L$ is reducible, contradicting the assumption that $L$ is irreducible. Therefore $a_j-a_i\notin S\left(b_i,b_j\right)$ for any $a_i$ and $a_j$ with $1\leq i\neq j\leq k$.
\end{proof}
\begin{remark}
For the Yangian $Y(\mathfrak{gl}_{l+1})$, Theorem 1.1 in \cite{Mo2} gives a necessary and sufficient condition for the tensor product of two evaluation modules to be irreducible. Moreover, Theorem 4.9 in \cite{NaTa} states that a tensor product of elementary $Y(\mathfrak{gl}_{l+1})$-modules is irreducible if and only if the tensor product of any two of those elementary modules is irreducible. These two theorems can be combined to obtain a more general result than our Theorem 5.21.
\end{remark}
In what follows, we prove that the local Weyl module $W(\pi)$ is isomorphic to an ordered tensor product of fundamental representations of $\yg$.
\begin{proposition}\label{wmiatpsp}
Let $\pi=\big(\pi_1(u),\ldots, \pi_{l}(u)\big)$, where $\pi_i\left(u\right)=\prod\limits_{j=1}^{m_i}\left(u-a_{i,j}\right)$. Let $S$ be the multiset of the roots of these polynomials. Let $a_1=a_{i,j}$ be one of the numbers in $S$ with the maximal real part and let $b_1=i$. Inductively, let $a_r=a_{s,t}\left(r\geq 2\right)$ be one of the numbers in $S\setminus\{a_1, \ldots, a_{r-1}\}\ (r\geq 2)$ with the maximal real part and $b_r=s$. Let $k=m_1+\ldots+m_l$. Then the ordered tensor product $L=V_{a_1}(\omega_{b_1})\otimes V_{a_2}(\omega_{b_2})\otimes\ldots\otimes V_{a_k}(\omega_{b_k})$ is a highest weight representation of $\yg$, and its associated polynomial is $\pi$.
\end{proposition}
\begin{proof}
Note that if $\operatorname{Re}(a_i)\geq \operatorname{Re}(a_j)$ for $1\leq i<j\leq k$, then $a_j-a_i\notin S(b_i,b_j)$. The rest of the proof follows from Theorem \ref{mtoyspl} and then Proposition \ref{vtv'hwv}.
\end{proof}

\begin{theorem}\label{wmiatpsp24}
The local Weyl module $W(\pi)$ of $\yg$ is isomorphic to the ordered tensor product $L$ as in Proposition \ref{wmiatpsp}.
\end{theorem}
\begin{proof}
On the one hand, $\operatorname{Dim}\big(W(\pi)\big)\leq \operatorname{Dim}\big(W(\lambda)\big)$ by Theorem \ref{ubodowm}; on the other hand, the dimension of $W(\pi)$ is $\geq \operatorname{Dim}\left(L\right)$ since $L$ is a quotient of $W(\pi)$. By Corollary \ref{dkrvawocsp}, $\operatorname{Dim}\big(W(\omega_i)\big)=\operatorname{Dim}V_{a}\left(\omega_i\right)$. Thus we have $\operatorname{Dim}\big(W(\lambda)\big)=\operatorname{Dim}\left(L\right)$, which implies that $\operatorname{Dim}\big(W(\pi)\big)=\operatorname{Dim}\left(L\right)$. Therefore $W(\pi)\cong L$.
\end{proof}

We conclude this section with the following remarks.  The methodology used in this paper can be also applied to describe the local Weyl modules and cyclicity condition for the tensor product of fundamental representations of $\yg$ when $\g$ is an exceptional simple Lie algebra over $\C$. The main difficulty in working with such Lie algebras is in finding the eigenvalues of $h_{r_j,k}$, $1\leq k\leq m_j$, on $v_{\sigma_j(\omega_i)}$. When $\g$ is a classical simple Lie algebra, $m_j\in\{0,1,2\}$, but $m_j$ can be greater than 2 when $\g$ is an exceptional simple Lie algebra. For instance, when $\g$ is of type $F_4$, $m_j\in\{0,1,2,3,4\}$ and when $\g$ is of type $E_8$, $m_j\in\{0,1,2,3,4,5,6\}$.
When $m_j\geq 3$, using defining relations of $\yg$ to compute the eigenvalues is more complicated.

It is possible to define modules for the Yangian $Y(\mathfrak{g})$ that would be similar to the global Weyl modules for loop and quantum loop algebras. It is normal to expect that they would share analogous properties to their classical and quantum loop counterparts, but this question deserves a separate publication. In particular, since results for global Weyl modules for loop algebras are not exactly the same as for quantum loop algebras, it is not clear to which of these two cases the Yangian case would be more similar.

A good number of papers have been written about the representation theory of twisted Yangians and the classification theorem for their finite dimensional representations is similar to the classification theorem for $Y(\mathfrak{g})$. It is possible to define for them modules that are analogous to the local Weyl modules for the Yangian $Y(\mathfrak{g})$. However, it is not clear if results similar to our main theorems could be obtained. Twisted Yangians are coideal subalgebras of the Yangian of $\mathfrak{gl}_n$ (not Hopf algebras) and the latter should be viewed using the RTT-presentation. It is not clear if this set of generators for $Y(\mathfrak{gl}_n)$ is well adapted to the study of local Weyl modules and the same could be said about the generators for twisted Yangians.
\section*{Acknowledgement}
The second author was supported by an NSERC Discovery Grant. The first author would like to thank the second author for supporting him through his NSERC Discovery Grant.

\end{document}